\documentclass[12pt]{article}
\usepackage{tikz}
\usetikzlibrary{intersections,calc,arrows}
\usepackage{amssymb}
\date{}
\makeatletter
\newcommand{\figcaption}[1]{\def\@captype{figure}\caption{#1}}
\newcommand{\tblcaption}[1]{\def\@captype{table}\caption{#1}}

\@addtoreset{equation}{section}
\makeatother
\newcommand{\qed}{\hbox{\rule[-2pt]{3pt}{6pt}}}

\begin{document}
\title {\bf Global bifurcation curves of nonlocal elliptic equations 
with oscillatory nonlinear term}

\author{{\bf Tetsutaro Shibata}
\\
{\small Laboratory of Mathematics, 
Graduate School of Advanced Science and Engineering,} 
\\
{\small Hiroshima University, Higashi-Hiroshima, 739-8527, Japan}
}

\maketitle
\footnote[0]{E-mail: tshibata@hiroshima-u.ac.jp}

\vspace{-0.5cm}

\begin{abstract}
We study the one-dimensional nonlocal elliptic equation of Kirchhoff type with 
oscillatory nonlinear term. 
We establish the precise asymptotic formulas for the bifurcation curves $\lambda(\alpha)$ 
as $\alpha \to \infty$ and $\alpha \to 0$, where $\alpha:= \Vert u_\lambda\Vert_\infty$ 
and $u_\lambda$ is the solution associated with 
$\lambda$.  We show that the second term of $\lambda(\alpha)$ is oscillatory 
as $\alpha \to \infty$.

\end{abstract}


\vspace{0.5cm}


\section{Introduction} 		      

We consider the following one-dimensional nonlocal elliptic equation

\begin{equation}
\left\{
\begin{array}{l}
-(b\Vert u'\Vert_2^2 + 1) u''(x)= \lambda(u(x)^p +  u(x)\sin^2 u(x)), \enskip 
x \in I:= (0,1),
\vspace{0.1cm}
\\
u(x) > 0, \enskip x\in I, 
\vspace{0.1cm}
\\
u(0) = u(1) = 0,
\end{array}
\right.
\end{equation}
where $p > 1, b\ge 0$ are given constants, $\lambda > 0$ is a bifurcation parameter and 
$\Vert \cdot\Vert_2$ denotes the usual $L^2$-norm.   

The purpose of  this paper is to establish the asymptotic formulas for bifurcation curves 
$\lambda = \lambda(\alpha)$ of 
(1.1) as $\alpha \to \infty$ 
to understand well how the oscillatory term gives effect to the bifurcation curves. 
Here $\alpha:= \Vert u_\lambda\Vert_\infty$ and $u_\lambda$ is a solution of (1.1) associated with 
$\lambda > 0$. 
When we consider the case where $b = 0$, we use the following 
notation to avoid the confusion:
\begin{equation}
\left\{
\begin{array}{l}
-v''(x)= \mu (v(x)^p + v(x)\sin^2 v(x)), \enskip 
x \in I,
\vspace{0.1cm}
\\
v(x) > 0, \enskip x\in I.
\vspace{0.1cm}
\\
v(0) = v(1) = 0,
\end{array}
\right.
\end{equation}
where $\mu > 0$ is the bifurcation parameter. It is well known by [12] that, for any given $\alpha > 0$, 
there exists a unique solution pair 
$(\mu, v_\alpha) \in \mathbb{R}_+ \times C^2(\bar{I})$ of (1.2) with 
$\alpha= \Vert v_\alpha\Vert_\infty$. 
Besides,  
$\mu$ is parameterized by $\alpha$ and 
a continuous function of $\alpha$ (cf. [12, Theorem 2.1]).
So we write $\mu = \mu(\alpha)$.

\noindent
Equation (1.1) is the nonlocal elliptic problem of Kirchhoff type motivated by the 
problem in [7]:
\begin{equation}
\left\{
\begin{array}{l}
-A\left(\displaystyle{\int_0^1} \vert u'(x)\vert^q dx\right)u''(x)= 
\lambda f(u(x)), \enskip 
x \in I,
\vspace{0.1cm}
\\
u(0) = u'(1) = 0,
\end{array}
\right.
\end{equation}
where $A = A(y)$, which is called Kirchhoff function (cf. [10, 15]), 
is a continuous function of $y \ge 0$. 
Nonlocal problems have been investigated by many authors and there are quite many manuscripts 
which treated the problems with the backgrounds in physics, 
biology, engineering and so on. We refer to [1--4, 6--9, 11, 13, 14], and the references therein. One of the main interests there are existence, nonexistence and the number of positive and 
nodal solutions. However, there seems to be a few works which considered (1.3) from 
a view-point of bifurcation problems. We refer to [16--21] and the references therein. 
As far as the author knows, there are no works which treat the nonlinear 
oscillatory 
eigenvalue problem such as (1.2).  
Therefore, there seems no works which treat nonlocal bifurcation problems with oscillatory nonlinear term, so our results here seem to be novel. Our approach are mainly 
the time-map method and the complicated calculation of definite integrals. 

The relationship between $\lambda(\alpha)$ and $\mu(\alpha)$ is as follows.  
Let $\alpha > 0$ be an arbitrary  given constant.
Assume that there exists a solution pair 
$(\lambda(\alpha), u_\alpha) \in \mathbb{R} \times C^2(\bar{I})$ with
$\Vert u_\alpha\Vert_\infty = \alpha$. Then we have 
\begin{eqnarray}
-u_\alpha''(x) 
= \frac{\lambda(\alpha)}{b\Vert u_\alpha'\Vert_2^2 + 1 }
(u_\alpha(x)^p + u_\alpha(x)\sin^2 u_\alpha(x)).
\end{eqnarray}
We note that $\Vert u_\alpha \Vert_\infty = \alpha$. Then we find that 
$u_\alpha = v_\alpha$ and 
$\frac{\lambda(\alpha)}{b\Vert u_\alpha'\Vert_2^2 + 1 } = \mu(\alpha)$, since the solution pair 
$(\mu(\alpha), v_\alpha) \in \mathbb{R}_+ \times C^2(\bar{I})$ of (1.2) 
with $\Vert v_\alpha \Vert_\infty = \alpha$ is unique (cf. [12]). This implies 
\begin{eqnarray}
\lambda(\alpha) =  (b\Vert v_\alpha'\Vert^2 + 1)\mu(\alpha).
\end{eqnarray}
Therefore, to obtain $\lambda(\alpha)$, we need to obtain both $\mu(\alpha)$ and 
$\Vert v_\alpha'\Vert_2$.

Now we state our results. We first consider the case $p > 2$.

\vspace{0.2cm}

\noindent
{\bf Theorem 1.1.} {\it  Consider (1.2). Let $p > 2$. Then as $\alpha \to \infty$,
\begin{eqnarray}
\mu(\alpha) &=& 2(p+1)\alpha^{1-p}
\left\{C_{0,p} + \left( C_1 + \frac12C_{11}\right)\alpha^{1-p} 
\right.
\\
&&
\left.\qquad 
+ \frac12(C_{12} + C_{21})\alpha^{-p}
+ \frac12C_{22}\alpha^{-(p+1)} + (C_2 + C_3)\alpha^{2(1-p)} + o(\alpha^{2(1-p)})\right\}^2,
\nonumber
\end{eqnarray}
where 
\begin{eqnarray}
C_{0,p}&:=& \int_0^1 \frac{1}{\sqrt{1-s^{p+1}}}ds, 
\\
C_1 &:=& -\frac{p+1}{8}\int_0^1 \frac{1-s^2}{(1-s^{p+1})^{3/2}}ds,
\\
C_{11} &:=& \frac{2}{p+1}\int_0^{\pi/2}\cos(2\alpha\sin^{2/(p+1}\theta)\sin^{(3-p)/(p+1)}
d\theta,
\\
C_{12} &:=& \frac{p-1}{2(p+1)}\int_0^{\pi/2}
(\sin2\alpha-\sin(2\alpha\sin^{2/(p+1)}\theta))\sin^{(1-p)/(p+1)}\theta d\theta
\\
&&\mbox{}+ \frac{p+1}{4}\int_0^1 \frac{1-s}{(1-s^{p+1})^{3/2}}\sin(2\alpha s)ds,
\nonumber
\\
C_{21} &:=& -\frac{1}{p+1}\int_0^{\pi/2} \sin(2\alpha\sin^{2/(p+1)}\theta)\sin^{(3-p)/(p+1)}\theta
d\theta 
\\
C_{22} &:=& \frac{4(p-1)}{p+1}\int_0^{\pi/2}(\cos2\alpha - \cos(2\alpha\sin^{2/(p+1)}\theta))
\sin^{(1-p)/(p+1)}\theta d\theta,
\\
C_2 &=& \frac{3(p+1)^2}{128}\int_0^1 \frac{(1-s^2)^2}{(1-s^{p+1})^{5/2}}ds,
\\
C_3 &=&  -\frac{3}{32}(p+1)^2\int_0^1 \left(\int_0^s \frac{1-y^2}{(1-y^{p+1})^{5/2}}dy\right)
\cos(2\alpha s)ds.
\end{eqnarray}
}
\vspace{0.2cm}

 \noindent
 {\bf Theorem 1.2.} {\it Consider (1.2). Let $p > 2$ and $\alpha \gg 1$. 
 Then the following asymptotic formula for $\Vert v_\alpha'\Vert_2^2$ holds. 
\begin{eqnarray}
 \Vert v_\alpha'\Vert_2^2 &=&  4\alpha^2\{G_0 + G_1\alpha^{1-p} + G_2\alpha^{-p} + G_3\alpha^{-(p+1)} 
+ G_4\alpha^{2(1-p)} + o(\alpha^{2(1-p)})\},
 \end{eqnarray}
 where 
\begin{eqnarray}
G_0 &:=& C_{0,p}E_{0,p}, 
\\
G_1 &:=& C_{0,p}E_1 +  \left( C_1 + \frac12C_{11}\right)E_{0,p},
\\
G_2 &:=& \frac12\left(C_{12} + C_{12}\right)E_{0,p} + C_{0,p}E_2,
\\
G_3 &:=&  \frac12C_{22}E_{0,p} + C_{0,p}E_3,
\\
G_4 &:=& (C_2 + C_3)E_{0,p} + C_{0,p}E_4 + \left( C_1 + \frac12C_{11}\right)E_1,
\\
E_{0,p} &:=& \int_0^1 \sqrt{1-s^{p+1}}ds, 
\\
E_1 &:=& \frac{p+1}{8}\int_0^1 \frac{1-s^4}{\sqrt{1-s^{p+1}}}ds,
\\
E_2 &:=& -\frac{1}{4}\int_0^{\pi/2}\{\sin2\alpha - \sin^{2/(p+1)}\theta\sin(2\alpha\sin^{2/(p+1)}\theta)
\}\sin^{(1-p)/(p+1)}\theta d\theta,
\\
E_3 &:=& -\frac{1}{8}\int_0^1 \{\cos2\alpha - \cos(2\alpha\sin^{2/(p+1)}\theta)
\}\sin^{(1-p)/(p+1)}d\theta,
\\
E_4 &:=& -\frac{(p+1)^2}{128}\int_0^1 \frac{(1-s^2)^2}{(1-s^{p+1})^{3/2}}ds,
\\
E_5 &:=& \frac{2}{p+1}\int_0^1 \frac{1-s^{p+1}}{\sqrt{1-s^4}}ds.
\end{eqnarray}
}

\vspace{0.2cm}

\noindent
{\bf Remark 1.3.} We should note that the order of the lower terms of $\mu(\alpha)$ in (1.6) 
changes according to $p$. Indeed, if we expand the bracket of the r.h.s. of (1.6), then the terms with 
$$
C_{0,p}^2, \alpha^{1-p}, \alpha^p, \alpha^{-(p+1)}, \alpha^{2(1-p)}, \alpha^{1-2p}
$$ 
appear.  Then for $\alpha \gg 1$, clearly, the first term is $C_{0,p}^2$ and the second is $\alpha^{1-p}$. Besides, we have 

\begin{equation}
\left\{
\begin{array}{l}
\alpha^{2(1-p)} \gg \alpha^{-p} \gg \alpha^{1-2p} \gg \alpha^{-(p+1)}  \qquad (1 < p < 2),
\vspace{0.1cm}
\\
\alpha^{-p} \sim \alpha^{2(1-p)} \gg \alpha^{-(p+1)} \sim \alpha^{1-2p} \qquad (p = 2),
\vspace{0.1cm}
\\
\alpha^{-p} \gg \alpha^{2(1-p)} \gg 
\alpha^{-(p+1)} \gg \alpha^{1-2p} \qquad (2 < p < 3),
\vspace{0.1cm}
\\
\alpha^{-p} \gg \alpha^{-(p+1)} \sim\alpha^{2(1-p)} \gg \alpha^{1-2p} \qquad (p = 3),
\vspace{0.1cm}
\\
 \alpha^{-p} \gg \alpha^{-(p+1)} \gg \alpha^{2(1-p)} \gg \alpha^{1-2p} \qquad (p > 3).
\end{array}
\right.
\end{equation}
Therefore, if $p > 2$, then the third term in the bracket of the r.h.s. of (1.6) is $\alpha^{-p}$. 
However, if $1 < p < 2$, then the third term is $\alpha^{2(1-p)}$. Moreover, if $p$ is very close to $1$, 
then $1-p \doteqdot 0$. Therefore, we have the sequence of the lower term, which are greater than 
$\alpha^{-p}$ in (1.6). In principle, it is possible to calculate them precisely. However, since the calculation is long and tedious, we do not carry out here. 

\vspace{0.2cm}

 \noindent
{\bf Theorem 1.4.} {\it Consider (1.2). 

\noindent
(i) Let $1 < p < 2$. Then as $\alpha \to \infty$,
\begin{eqnarray}
\mu(\alpha) = 2(p+1)\alpha^{1-p}
\left\{C_{0,p} + \left( C_1 + \frac12C_{11}\right)\alpha^{1-p} + (C_2 + C_3)\alpha^{2(1-p)} + o(\alpha^{2(1-p)})\right\}^2.
\end{eqnarray}
(ii) Let $p = 2$. Then as $\alpha \to \infty$,
\begin{eqnarray}
\mu(\alpha) &=& 6\alpha^{-1}
\\
&&\times
\left\{C_{0,p} + \left( C_1 + \frac12C_{11}\right)\alpha^{-1} 
+ \left(\frac12C_{12} + \frac12C_{21} + C_2 + C_3\right)\alpha^{-2} + o(\alpha^{2(1-p)})\right\}^2.
\nonumber
\end{eqnarray}

}

\vspace{0.2cm}

 \noindent
{\bf Theorem 1.5.} {\it Consider (1.2). 

\noindent
(i) Let $1 < p < 2$. Then as $\alpha \to \infty$,
\begin{eqnarray}
 \Vert v_\alpha'\Vert_2^2 &=&  4\alpha^2\{G_0 + G_1\alpha^{1-p} + G_4\alpha^{2(1-p)} 
 + G_2\alpha^{-p} + o(\alpha^{2(1-p)})\}.
 \end{eqnarray}
 (ii) Let $p = 2$. Then as $\alpha \to \infty$,
\begin{eqnarray}
 \Vert v_\alpha'\Vert_2^2 &=&  4\alpha^2\{G_0 + G_1\alpha^{-1} + (G_2 + G_4)\alpha^{-2} 
 + o(\alpha^{-2})\}.
 \end{eqnarray}
}
 
 \vspace{0.2cm}

Theorems 1.4 and 1.5 are obtained directly from Theorems 1,1 and 1.2.  So we omit the proofs. 
 
 We now consider (1.1).

\vspace{0.2cm}
 
 \noindent
 {\bf Theorem 1.6.} {\it Consider (1.1) with $b > 0$. 
 
 \noindent
 (i) Let $p > 2$ and $\alpha \gg 1$. 
 Then the following asymptotic formula for $\lambda(\alpha)$ holds. 
\begin{eqnarray}
\lambda(\alpha) &=&  2(p+1)\alpha^{1-p}
\left\{C_{0,p} + \left( C_1 + \frac12C_{11}\right)\alpha^{1-p} 
\right.
\\
&&
\left.\qquad \qquad \qquad 
+ \frac12(C_{12} + C_{21})\alpha^{-p}
+ \frac12C_{22}\alpha^{-(p+1)} + C_2\alpha^{2(1-p)} + o(\alpha^{2(1-p)})\right\}^2
\nonumber
\\
&&\times
\left\{4b\alpha^2\{G_0 + G_1\alpha^{1-p} + G_2\alpha^{-p} + G_3\alpha^{-(p+1)} 
+ G_4\alpha^{2(1-p)} + o(\alpha^{2(1-p)})\} + 1\right\}.
\nonumber
\end{eqnarray}
(ii) Let $p = 2$. Then as $\alpha \to \infty$, 
\begin{eqnarray}
\lambda(\alpha) &=&  6\alpha^{-1}
\left\{C_{0,p} + \left( C_1 + \frac12C_{11}\right)\alpha^{-1} 
+ \left(\frac12C_{12} + \frac12C_{21} + C_2 + C_3\right)\alpha^{-2} + o(\alpha^{2(1-p)})\right\}^2
\nonumber
\\
&&
\times \left\{4b\alpha^2\{G_0 + G_1\alpha^{-1} + (G_2 + G_4)\alpha^{-2} 
 + o(\alpha^{-2})\} + 1\right\}.
\end{eqnarray}

\noindent
 (iii) Let $1 < p < 2$ Then as $\alpha \to \infty$, 
 \begin{eqnarray}
\lambda(\alpha) &=&  2(p+1)\alpha^{1-p}
\left\{C_{0,p} + \left( C_1 + \frac12C_{11}\right)\alpha^{1-p} 
+ C_2\alpha^{2(1-p)} + o(\alpha^{2(1-p)})\right\}^2
\nonumber
\\
&&\times
\left\{4b\alpha^2\{G_0 + G_1\alpha^{1-p} + G_4\alpha^{2(1-p)} + o(\alpha^{2(1-p)})\} + 1\right\}.
\end{eqnarray}
 } 
 
\vspace{0.2cm}

We see from Theorem 1.6 that, roughly speaking, the asymptotic behaviors of $\lambda(\alpha)$ as $\alpha \to \infty$ are:
\begin{eqnarray}
\lambda(\alpha) \sim \alpha^{3-p}.
\end{eqnarray}

We obtain Theorem 1.6 immediately by (1.5), 
Theorems 1.1, 1.2, 1.4 and 1.5. So we omit the proof. 

Now we establish the asymptotic formulas for $\mu(\alpha)$ as $\alpha \to 0$ to understand 
the entire structure of $\mu(\alpha)$. We put 

\begin{eqnarray}
H_2&:=& -\frac{2}{p+1}
\int_0^1 \frac{1-s^{p+1}}{(1-s^4)^{3/2}}ds, 
\\
H_n &:=& -2^{2n-2}(-1)^n\left\{\frac{1}{(2n-1)!}\int_0^1 \frac{1-s^{2n-1}}{(1-s^4)^{3/2}}ds 
-\frac{1}{(2n)!}\int_0^1 \frac{1-s^{2n}}{(1-s^4)^{3/2}}ds\right\} 
\end{eqnarray}
for $n \ge 3$. Furthermore, let
\begin{eqnarray}
L_1 &:=& -\frac{p+1}{8}\int_0^1 
\frac{1-s^4}{(1-s^{p+1})^{3/2}} ds,
\\
L_2&:=& -\frac12\int_0^1 \frac{1}{\sqrt{1-s^{p+1}}}K(s)ds,
\\
K(s) &:=& - 2^3(p+1)\left\{\frac{1}{5!}\frac{1-s^5}{1-s^{p+1}} 
- \frac{1}{6!}\frac{1-s^6}{1-s^{p+1}} 
+ O(\alpha^{7-p})\right\}.
\end{eqnarray}

\vspace{0.2cm}

 \noindent
{\bf Theorem 1.7.} {\it Consider (1.2). 

\noindent
(i) Let $1 < p < 3$. Then as $\alpha \to 0$,
\begin{eqnarray}
\mu(\alpha) &=& 2(p+1)\alpha^{1-p}\left\{C_{0,p} + L_1\alpha^{3-p} 
+ L_2\alpha^{5-p} + O(\alpha^{7-p})\right\}^2.
\end{eqnarray}

\noindent
(ii) Let $p = 3$. Then as $\alpha \to 0$,
\begin{eqnarray}
\mu(\alpha) = 4\alpha^{-2}\left\{C_{0,3} +\frac12H_3\alpha^2 + O(\alpha^4)\right\}^2.
\end{eqnarray}

\noindent
(iii) Let $3 < p \le 5$. Then as $\alpha \to 0$, 
\begin{eqnarray}
\mu(\alpha)
&=& 8\alpha^{-2}\left\{C_{0,3} + H_2\alpha^{p-3} + H_3\alpha^2 + O(\alpha^4)\right\}^2. 
\end{eqnarray}

\noindent
(iv) Assume that $p > 5$. Then as $\alpha \to 0$,
\begin{eqnarray}
\mu(\alpha) &=& 8\alpha^{-2}\left\{C_{0,3} + H_3\alpha^{2} + o(\alpha^{2})
\right\}^2.
\end{eqnarray}

}
\vspace{0.2cm}

Finally, we establish the asymptotic formulas for $\lambda(\alpha)$ as $\alpha \to 0$. 

\vspace{0.2cm}

 \noindent
{\bf Theorem 1.8.} {\it Consider (1.1). 

\noindent
(i) Let $1 < p < 3$. Then as $\alpha \to 0$, 
 \begin{eqnarray}
 \lambda(\alpha) &=& 2(p+1)\alpha^{1-p}\left\{C_{0,p} + L_1\alpha^{3-p} 
+ L_2\alpha^{5-p} + O(\alpha^{7-p})\right\}^2
\\
&&\times
\left\{4b\alpha^2
\left\{E_{0,p}C_{0,p} + (E_{0,p}L_1+C_{0,p}E_1)\alpha^{3-p} + o(\alpha^{3-p})\right\} + 1\right\}.
\nonumber
 \end{eqnarray}
  (ii) Let $p = 3$. Then as $\alpha \to 0$, 
  \begin{eqnarray}
 \lambda(\alpha) &=& 4\alpha^{-2}(1 + 4bE_{0.3}C_{0,3}\alpha^2 + o(\alpha^2))
 \left\{C_{0,3} + \frac12H_3\alpha^2 + O(\alpha^4)\right\}^2.
 \end{eqnarray}
 (iii) Let $3 < p \le 5$. Then as $\alpha \to 0$, 
\begin{eqnarray}
\lambda(\alpha) &=& 8\alpha^{-2}
\left\{C_{0,3} + H_2\alpha^{p-3} + H_3\alpha^2 + O(\alpha^4)\right\}^2
\\
&&\times
\left[4b\alpha^{2}\left\{C_{0,3} + H_2\alpha^{p-3} + H_3\alpha^2 + O(\alpha^4)\right\}
\left\{E_{0,3} + E_5\alpha^{p-3}(1 + o(1))\right\} + 1\right].
\nonumber
\end{eqnarray}
(iv) Let $p > 5$. Then as $\alpha \to 0$, 
\begin{eqnarray}
\lambda(\alpha) &=& 8\alpha^{-2}\left\{C_{0,3} + H_3\alpha^2 + o(\alpha^2)\right\}^2
\\
&&\times
\left[4b\alpha^{2}\left\{C_{0,3} + H_3\alpha^2 + o(\alpha^2)\right\}
\left\{E_{0,3} + E_5\alpha^{p-3}(1 + o(1))\right\} + 1\right].
\nonumber
\end{eqnarray}
}

\vspace{0.2cm}
By Theorem 1.8, we see that as $\alpha \to 0$, 

\begin{equation}
\lambda(\alpha) \sim \left\{
\begin{array}{l}
\alpha^{1-p} \enskip (1 < p \le 3),
\vspace{0.1cm}
\\
\alpha^{-2} \enskip (p > 3).
\end{array}
\right.
\end{equation}

\section{Proofs of Theorems 1.1 and 1.2}

In this section, let $p > 2$ and we consider (1.2). In what follows, $C$ denotes various positive 
constants independent of $\alpha \gg 1$. 
By [5], we know that if $v_\alpha$ is a solution of (1.2), then $v_\alpha$ satisfies 
\begin{eqnarray}
v_\alpha(x) &=& v_\alpha(1-x), \quad 0 \le x\le \frac12,
\\
\alpha&:=& \Vert v_\alpha\Vert_\infty = v_\alpha\left(\frac12\right),
\\
v_\alpha'(x) &>& 0, \quad 0 \le x < \frac12.
\end{eqnarray}
We put 
\begin{eqnarray}
f(\theta) &:=& \theta^p + \theta\sin^2 \theta, 
\\
F(\theta) &:=& \int_0^\theta f(y)dy = \frac{1}{p+1}\theta^{p+1} + 
\frac14\theta^2 - \frac14\theta\sin 2\theta - \frac18 \cos 2\theta + \frac18.
\end{eqnarray}
Let $\alpha > 0$ be an arbitrary given constant. We write $\mu = \mu(\alpha)$ and $v_\alpha:= v_{\mu(\alpha)}$ in what follows. By (1.2), for $x \in \bar{I}$, we have 
\begin{eqnarray}
\{v_\alpha''(x) + \mu f(v_\alpha(x)\}v_\alpha'(x) = 0.
\end{eqnarray}
By this and (2.2),  for $x \in \bar{I}$, we have 
\begin{eqnarray}
\frac12v_\alpha'(x)^2 + \mu F(v_\alpha(x)) 
= \mbox{constant} = \mu F\left(v_\alpha\left(\frac{1}{2}\right)\right) = \mu F(\alpha).
\end{eqnarray}
By this and (2.3), for $0 \le x \le 1/2$, we have 
\begin{eqnarray}
v_\alpha'(x) &=& \sqrt{2\mu(F(\alpha)-F(v_\alpha(x)))} 
\\
&=&
\sqrt{\frac{2\mu}{p+1}}
\sqrt{(\alpha^{p+1} - v_\alpha(x)^{p+1}) + \frac{p+1}{4}(\alpha^2- v_\alpha(x)^2) - 
A_\alpha(v_\alpha(x))  - B_\alpha(v_\alpha(x))},
\nonumber
\end{eqnarray}
where 
\begin{eqnarray}
A_\alpha(v_\alpha(x))&:=& \frac{p+1}{4}(\alpha\sin 2\alpha - v_\alpha(x)
\sin(2v_\alpha(x))), 
\\
B_\alpha(v_\alpha(x))&:=& \frac{p+1}{8}(\cos 2\alpha - \cos (2v_\alpha(x))).
\end{eqnarray}
Note that $A_\alpha(v_\alpha(x)) \ll \alpha^2, B_\alpha(v_\alpha(x)) \ll \alpha^2$. 
By this and putting $v_\alpha(x) = \alpha s$, we have
\begin{eqnarray}
\frac12&=& \int_0^{1/2} 1 dx 
\\
&=& \sqrt{\frac{p+1}{2\mu}}\int_0^{1/2} \frac{v_\alpha'(x)dx}
{ \sqrt{(\alpha^{p+1} - v_\alpha(x)^p) + \frac{p+1}{4}(\alpha^2- v_\alpha(x)^2) - 
A_\alpha(v_\alpha(x))  - B_\alpha(v_\alpha(x))}
}
\nonumber
\\
&=& \sqrt{\frac{p+1}{2\mu}}\alpha^{(1-p)/2}\int_0^1 \frac{ds}{\sqrt{(1-s^{p+1}) +\frac{p+1}{4}\alpha^{1-p}(1-s^2)- \frac{1}{\alpha^{p+1}}A_\alpha(\alpha s) - \frac{1}{\alpha^{p+1}}B_\alpha(\alpha s)}}
\nonumber
\\
&=& 
\sqrt{\frac{p+1}{2\mu}}\alpha^{(1-p)/2}\int_0^1 \frac{1}{\sqrt{1-s^{p+1}}}
\frac{ds}
{\sqrt{1  + \frac{p+1}{4}\alpha^{1-p}\frac{1-s^2}{1-s^{p+1}} 
-\frac{1}{\alpha^{p+1}}\frac{A_\alpha(\alpha s)}{1-s^{p+1}} 
- \frac{1}{\alpha^{p+1}} \frac{B_\alpha(\alpha s)}{1-s^{p+1}}}
\nonumber
}.
\end{eqnarray} 
This along with Taylor expansion implies that 
\begin{eqnarray}
\sqrt{\mu} &=& \sqrt{2(p+1)}\alpha^{(1-p)/2}
\\
&&\times\int_0^1 \frac{1}{\sqrt{1-s^{p+1}}}
\left\{1 - \frac{p+1}{8}\alpha^{1-p}\frac{1-s^2}{1-s^{p+1}} 
+ \frac12 \frac{1}{\alpha^{p+1}}\frac{A_\alpha(\alpha s)}{1-s^{p+1}} + \frac12
\frac{1}{\alpha^{p+1}} \frac{B_\alpha(\alpha s)}{1-s^{p+1}}
\right.
\nonumber
\\
&& \left.+ \frac{3}{8}\left(\frac{p+1}{4}\alpha^{1-p}\frac{1-s^2}{1-s^{p+1}}\right)^2 
- \frac{3}{16}(p+1)\alpha^{-2p}\frac{1-s^2}{(1-s^{p+1})^2}A_\alpha(\alpha s) + o(\alpha^{2(1-p)})
\right\}ds
\nonumber
\\
&=&  \sqrt{2(p+1)}\alpha^{(1-p)/2}\left[C_{0,p} + C_1\alpha^{1-p} + I + II + C_2\alpha^{2(1-p)} 
+ III + o(\alpha^{2(1-p)})\right],
\nonumber
\end{eqnarray}
where 
\begin{eqnarray}
I &=& \frac{1}{2}\alpha^{-(p+1)}I_1 := \frac{1}{2}\alpha^{-(p+1)}\int_0^1 
\frac{A_\alpha(\alpha s )}{(1-s^{p+1})^{3/2}}ds,
\\
II &=& \frac{1}{2}\alpha^{-(p+1)}II_1 := \frac{1}{2}\alpha^{-(p+1)}
\int_0^1 \frac{B_\alpha(\alpha s )}
{(1-s^{p+1})^{3/2}}ds,
\\
III &=& -\frac{3}{16}(p+1)\alpha^{-2p}\int_0^1 \frac{1-s^2}{(1-s^{p+1})^{5/2}}A_\alpha(\alpha s)ds.
\end{eqnarray}

\vspace{0.2cm}

\noindent
{\bf Lemma 2.1.} {\it Let $\alpha \gg 1$. Then 
\begin{eqnarray}
I_1&=& \int_0^1 \frac{A_\alpha(\alpha s)}{(1-s^{p+1})^{3/2}}ds
= C_{11}\alpha^2 + C_{12}\alpha,
\\
II_1 &=& \int_0^1 \frac{B_\alpha(\alpha s )}
{(1-s^{p+1})^{3/2}}ds = C_{21}\alpha + C_{22}.
\end{eqnarray}
}
\noindent
{\it Proof.} We first note that the definite integrals $C_{11}, C_{12}, C_{21}, C_{22}$ exist, since
we have $-1 < (1-p)/(p+1) < (3-p)/(p+1)$. We first prove (2.16). We put $s:= \sin^{2/(p+1)}\theta$. Then by integration by parts, 
we have 
\begin{eqnarray}
I_{1} &=& \frac{p+1}{4}\alpha\int_0^1 \frac{1}{\sqrt{1-s^{p+1}}}\frac{\sin2\alpha - \sin(2\alpha s)}
{1-s^{p+1}}ds 
\\
&& +  \frac{p+1}{4}\alpha\int_0^1 \frac{(1-s)}{(1-s^{p+1})^{3/2}}\sin(2\alpha s)ds
\nonumber
\\
&=& \frac{1}{2}\alpha\int_0^{\pi/2} \frac{1}{\cos^2\theta}
\left[\left\{\sin2\alpha - \sin(2\alpha\sin^{2/(p+1)}\theta)\right\}\sin^{(1-p)/(p+1)}
\theta\right]d\theta
\nonumber
\\
&&  + \frac{p+1}{4}\alpha\int_0^1 \frac{(1-s)}{(1-s^{p+1})^{3/2}}\sin(2\alpha s)ds
\nonumber
\\
&=& \frac{1}{2}\alpha\int_0^{\pi/2} (\tan\theta)'
\left[\left\{\sin2\alpha - \sin(2\alpha\sin^{2/(p+1)}\theta)\right\}\sin^{(1-p)/(p+1)}
\theta\right]d\theta
\nonumber
\\
&&  + \frac{p+1}{4}\alpha\int_0^1 \frac{(1-s)}{(1-s^{p+1})^{3/2}}\sin(2\alpha s)ds
\nonumber
\\
&=& \frac12\alpha\left[\tan\theta
\left[\left\{\sin2\alpha - \sin(2\alpha\sin^{2/(p+1)}\theta\right\}\sin^{(1-p)/(p+1)}
\theta\right]\right]^{\pi/2}_0 \cdots\cdots \cdots \cdots(*)
\nonumber
\\
&&- \frac{1}{2}\alpha\int_0^{\pi/2}\frac{\sin\theta}{\cos\theta} 
\left\{-\frac{4}{p+1}\alpha 
\cos(2\alpha\sin^{2/(p+1}\theta)\sin^{(2-2p)/(p+1)}\theta \cos\theta
\right.
\nonumber
\\
&& \qquad \qquad \qquad \left.
-\frac{p-1}{p+1}(\sin2\alpha-\sin(2\alpha\sin^{2/(p+1)}\theta))\sin^{-2p/(p+1)}\theta\cos\theta
\right\}d\theta
\nonumber
\\
&& +  \frac{p+1}{4}\alpha\int_0^1 \frac{(1-s)}{(1-s^{p+1})^{3/2}}\sin(2\alpha s)ds
\nonumber
\\
&=& \frac{2}{p+1}\alpha^2 \int_0^{\pi/2}\cos(2\alpha\sin^{2/(p+1}\theta)\sin^{(3-p)/(p+1)}
d\theta 
\nonumber
\\
&&+ \frac{p-1}{2(p+1)}\alpha\int_0^{\pi/2}
(\sin2\alpha-\sin(2\alpha\sin^{2/(p+1)}\theta))\sin^{(1-p)/(p+1)}\theta d\theta
\nonumber
\\
&& +  \frac{p+1}{4}\alpha\int_0^1 \frac{(1-s)}{(1-s^{p+1})^{3/2}}\sin(2\alpha s)ds
\nonumber
\\
&:=& C_{11}\alpha^2 + C_{12}\alpha.
\nonumber
\end{eqnarray} 
We remark that by l'H\^opital's rule and direct calculation, we easily obtain that $(*)$ in (2.18) 
and $(**)$ in 
(2.19) below are equal to $0$. Next, we put $s:= \sin^{2/(p+1)}\theta$. Then by integration by parts, 
we have 
\begin{eqnarray}
II_{1} &=& \frac{1}{4}\int_0^{\pi/2} \frac{1}{\cos^2\theta}
\left\{\cos2\alpha - \cos(2\alpha\sin^{2/(p+1)}\theta)\right\}\sin^{(1-p)/(p+1)}\theta d\theta
\\
&=&  \frac{1}{4}\int_0^{\pi/2} (\tan\theta)'
\left\{\cos2\alpha - \cos(2\alpha\sin^{2/(p+1)}\theta)\right\}\sin^{(1-p)/(p+1)}\theta d\theta
\nonumber
\\
&=& \frac14\left[ \tan\theta
\left\{\cos2\alpha - \cos(2\alpha\sin^{2/(p+1)}\theta)\right\}\sin^{(1-p)/(p+1)}\theta 
\right]_0^{\pi/2} \cdots\cdots \cdots \cdots(**)
\nonumber
\\
&& \enskip -\frac{1}{p+1}\alpha\int_0^{\pi/2} \sin(2\alpha\sin^{2/(p+1)}\theta)\sin^{(3-p)/(p+1)}\theta
d\theta 
\nonumber
\\
&& \enskip + \frac{4(p-1)}{p+1}\int_0^{\pi/2}(\cos2\alpha - \cos(2\alpha\sin^{2/(p+1)}\theta))
\sin^{(1-p)/(p+1)}\theta d\theta
\nonumber
\\
&=& C_{21}\alpha + C_{22}.
\nonumber
\end{eqnarray}
Thus the proof is complete. \qed

\vspace{0.2cm}

\noindent
{\bf Lemma 2.2.} {\it Let $\alpha \gg 1$. Then 
\begin{eqnarray}
III &=& C_3 \alpha^{2(1-p)} + o(\alpha^{2(1-p)}).
\end{eqnarray}
}
{\it Proof.} by (2.9) and (2.15), we have 
\begin{eqnarray}
III &=& -\frac{3}{64}(p+1)^2\alpha^{-2p}\int_0^1 \frac{1-s^2}{(1-s^{p+1})^{5/2}}
\left\{\alpha \sin 2\alpha - \alpha s \sin(2\alpha s)\right\}ds 
\\
& =&   -\frac{3}{64}(p+1)^2\alpha^{-2p+1}\int_0^1 \frac{1-s^2}{(1-s^{p+1})^{5/2}}
\left\{\sin 2\alpha - \sin (2\alpha s) \right\}ds 
\nonumber
\\
&&  -\frac{3}{64}(p+1)^2\alpha^{-2p+1}\int_0^1 \frac{(1-s^2)(1-s)}{(1-s^{p+1})^{5/2}}
\sin(2\alpha s)ds.
\nonumber
\\
& :=&   -\frac{3}{64}(p+1)^2\alpha^{-2p+1}III_1 + O(\alpha^{-2p+1}).
\nonumber
\end{eqnarray}
We show that $III_1 \sim \alpha$. 
We note that $(1-y^2) /(1-y^{p+1})^{5/2} \le (1-y^2)^{-3/2}$ for 
$0 \le y \le 1$. By this and integration by parts, we have 
\begin{eqnarray}
III_1 &=& 
\lim_{\epsilon \to 0}\int_0^{1-\epsilon} 
\frac{d}{ds}\left(\int_0^s\frac{1-y^2}{(1-y^{p+1})^{5/2}}dy\right)
\left\{\sin 2\alpha - \sin (2\alpha s) \right\}ds 
\\
&=& 
\lim_{\epsilon \to 0}\left[\left(\int_0^s \frac{1-y^2}{(1-y^{p+1})^{5/2}}dy\right)\{\sin 2\alpha 
- \sin(2\alpha s)\}\right]_0^{1-\epsilon} 
\nonumber
\\
&&+ 2\alpha\lim_{\epsilon \to 0}\int_0^{1-\epsilon} \left(\int_0^s \frac{1-y^2}{(1-y^{p+1})^{5/2}}dy\right)\cos(2\alpha s)ds
\nonumber
\\
&=& 2\alpha(1 + o(1))\int_0^1 \left(\int_0^s \frac{1-y^2}{(1-y^{p+1})^{5/2}}dy\right)\cos(2\alpha s)ds.
\nonumber
\end{eqnarray}
By this and (2.21), we have (2.20). Thus the proof is complete. \qed

\vspace{0.2cm}

\noindent
{\it Proof of Theorem 1.1.} 
By (2.12) and Lemma 2.1, for $\alpha \gg 1$, we obtain 
\begin{eqnarray}
\sqrt{\mu} &=& \sqrt{2(p+1)}\alpha^{(1-p)/2}
\left\{C_{0,p} + ( C_1 + \frac12C_{11})\alpha^{1-p} 
\right.
\\
&&
\left.\qquad 
+ \frac12(C_{12} + C_{21})\alpha^{-p}
+ \frac12C_{22}\alpha^{-(p+1)} + (C_2 + C_3)\alpha^{2(1-p)} + o(\alpha^{2(1-p)})\right\}.
\nonumber
\end{eqnarray}
By this, we obtain Theorem 1.1. Thus the proof is complete. \qed

\vspace{0.2cm}

We next prove Theorem 1.2. 

\vspace{0.2cm}

\noindent
{\bf Lemma 2.3.} {\it Let $v_\alpha$ be the solution of (1.2) associated with $\mu > 0$ such that 
$\Vert v_\alpha\Vert_\infty = \alpha > 0$. Then for $\alpha \gg 1$
\begin{eqnarray}
\Vert v_\alpha'\Vert_2^2 &=& 4\alpha^2\{G_0 + G_1\alpha^{1-p} + G_2\alpha^{-p} + G_3\alpha^{-(p+1)} 
+ G_4\alpha^{2(1-p)} + o(\alpha^{2(1-p)})\}.
\end{eqnarray}
}
\noindent
{\it Proof.} By (2.8), putting $v_\alpha(x) = \alpha s$ and Taylor expansion, we obtain
\begin{eqnarray}
\Vert v_\alpha'\Vert_2^2 &=& 2\int_0^{1/2} v_\alpha'(x)v_\alpha'(x)dx 
\\
&=& 2\sqrt{\frac{2\mu}{p+1}}
\nonumber
\\
&&\times\int_0^{1/2}
\sqrt{(\alpha^{p+1} - v_\alpha(x)^p) + \frac{p+1}{4}(\alpha^2- v_\alpha(x)^2) - 
A_\alpha(v_\alpha(x))  - B_\alpha(v_\alpha(x))}v_\alpha'(x)dx
\nonumber
\\
&&=  2\sqrt{\frac{2\mu}{p+1}}\alpha^{(p+3)/2}\int_0^{1}
\sqrt{1-s^{p+1}}
\nonumber
\\
&&\qquad \qquad \qquad \times 
\sqrt{1 +\frac{p+1}{4}\alpha^{1-p}\frac{1-s^2}{1-s^{p+1}}- 
\frac{1}{\alpha^{p+1}}\frac{A_\alpha(\alpha s)}{1-s^{p+1}} 
- \frac{1}{\alpha^{p+1}}\frac{B_\alpha(\alpha s)}{1-s^{p+1}}}ds
\nonumber
\\
&&= 2\sqrt{\frac{2\mu}{p+1}}\alpha^{(p+3)/2}\int_0^{1}
\sqrt{1-s^{p+1}}\left\{1 +\frac{p+1}{8}\alpha^{1-p}\frac{1-s^2}{1-s^{p+1}}- 
\frac{1}{2\alpha^{p+1}}\frac{A_\alpha(\alpha s)}{1-s^{p+1}} 
\right.
\nonumber
\\
&& \qquad \qquad 
- \frac{1}{2\alpha^{p+1}}\frac{B_\alpha(\alpha s)}{1-s^{p+1}} - \frac{(p+1)^2}{128}
\alpha^{2(1-p)}\left(\frac{1-s^2}{1-s^{p+1}}\right)^2 
\nonumber
\\
&& \left.\qquad \qquad 
+ \frac{1}{64}(p+1)^2\alpha^{-2p}\frac{1-s^2}{1-s^{p+1}}(\alpha\sin2\alpha 
- \alpha s\sin(2\alpha s))
+ o(\alpha^{2(1-p)})\right\}ds.
\nonumber
\end{eqnarray} 
By putting $s = \sin^{2/(p+1)}\theta$, we have
\begin{eqnarray}
\int_0^1 \frac{A_\alpha(\alpha s)}{\sqrt{1-s^{p+1}}}ds
&=& \frac{p+1}{4}\alpha\int_0^1 \frac{\sin2\alpha - s\sin(2\alpha s)}{\sqrt{1-s^{p+1}}}ds
\\
&=& \frac{1}{2}\alpha\int_0^{\pi/2}\{\sin2\alpha - \sin^{2/(p+1)}\theta\sin(2\alpha\sin^{2/(p+1)}\theta)
\}\sin^{(1-p)/(p+1)}\theta d\theta,
\nonumber 
\\
\int_0^1 \frac{B_\alpha(\alpha s)}{\sqrt{1-s^{p+1}}}ds
&=& \frac{p+1}{8}\int_0^1 \frac{\cos2\alpha - \cos(2\alpha s)}{\sqrt{1-s^{p+1}}}ds
\\
&=& \frac{1}{4}\int_0^1 \{\cos2\alpha - \cos(2\alpha\sin^{2/(p+1)}\theta)
\}\sin^{(1-p)/(p+1)}d\theta.
\nonumber
\end{eqnarray}
By (2.25)--(2.27), we have 
\begin{eqnarray}
\Vert v_\alpha'\Vert_2^2 &=& 2\sqrt{\frac{2\mu}{p+1}}\alpha^{(p+3)/2}
\left\{E_{0,p} + E_1\alpha^{1-p} + E_2 \alpha^{-p}  + E_3\alpha^{-(p+1)} 
\right.
\\
&&
\qquad \qquad \qquad \qquad \left.
+ E_4 \alpha^{2(1-p)} + o(\alpha^{2(1-p)})\right\}.
\nonumber
\end{eqnarray}
By this, (2.23)--(2.28), we have 
\begin{eqnarray}
\Vert v_\alpha'\Vert_2^2 &=& 4\alpha^2
\left\{C_{0,p} + \left( C_1 + \frac12C_{11}\right)\alpha^{1-p} 
+ \frac12(C_{12} + C_{21})\alpha^{-p}
\right.
\\
 &&\qquad \qquad \qquad\left.
+ \frac12C_{22}\alpha^{-(p+1)} + (C_2 + C_3)\alpha^{2(1-p)} + o(\alpha^{2(1-p)})\right\}
\nonumber
\\
&& \qquad \times\left\{E_{0,p} + E_1\alpha^{1-p} + E_2\alpha^{-p} + E_3\alpha^{-(p+1)} 
+ E_4\alpha^{2(1-p)} + o(\alpha^{2(1-p)})\right\}
\nonumber
\\
&=& 4\alpha^2\{G_0 + G_1\alpha^{1-p} + G_2\alpha^{-p} + G_3\alpha^{-(p+1)} 
+ G_4\alpha^{2(1-p)} + o(\alpha^{2(1-p)})\}.
\nonumber
\end{eqnarray}
This implies (2.27). Thus the proof is complete. \qed

\section{Proof of Theorem 1.7}

In this section, let $0 < \alpha \ll 1$. We put $w_\alpha:= v_\alpha/\alpha$. By (2.5) and Taylor expansion, we have 
\begin{eqnarray}
F(\alpha) &=& \frac{1}{p+1}\alpha^{p+1} + \frac14\alpha^2 - \frac14\alpha
\left\{2\alpha - \frac{1}{3!}(2\alpha)^3 + \sum_{n = 3}^\infty \frac{(-1)^{n-1}}
{(2n-1)!}(2\alpha)^{2n-1}\right\}
\\
&& - \frac{1}{8}\left\{1 - \frac{1}{2!}(2\alpha)^2 + \frac{1}{4!}(2\alpha)^4 
+ \sum_{n=3}^\infty \frac{(-1)^n}{(2n)!}(2\alpha)^{2n}\right\} + \frac18,
\nonumber
\\
F(v_\alpha) &=&  \frac{1}{p+1}\alpha^{p+1}w_\alpha^{p+1} + \frac14\alpha^2 w_\alpha(x)^2
- \frac14\alpha w_\alpha
\left\{2\alpha w_\alpha  - \frac{1}{3!}(2\alpha w_\alpha)^3 
+ \sum_{n = 3}^\infty \frac{(-1)^{n-1}}
{(2n-1)!}(2\alpha w_\alpha)^{2n-1}\right\}
\nonumber
\\
&& - \frac{1}{8}\left\{1 - \frac{1}{2!}(2\alpha w_\alpha)^2 + \frac{1}{4!}(2\alpha w_\alpha)^4 
+ \sum_{n=3}^\infty \frac{(-1)^n}{(2n)!}(2\alpha w_\alpha)^{2n}\right\} + \frac18.
\end{eqnarray}
By the same argument as that to obtain (2.8), for $0 \le x \le 1$, we have 
\begin{eqnarray}
\frac12\alpha^2 w_\alpha'(x)^2 &=& \mu
\left\{ \frac{1}{p+1}\alpha^{p+1}(1-w_\alpha(x)^{p+1}) + \frac14\alpha^4
(1-w_\alpha(x)^4) \right.
\\
&&\mbox{} + \frac14\alpha\sum_{n=3}^\infty
\frac{(-1)^n}{(2n-1)!}2^{2n-1}\alpha^{2n-1}(1-w_\alpha(x)^{2n-1}) 
\nonumber
\\
&&\mbox{}\left.
- \frac18\sum_{n=3}^\infty \frac{(-1)^n}{(2n)!}2^{2n}\alpha^{2n}(1 - w_\alpha(x)^{2n})\right\}.
\nonumber
\end{eqnarray}
We put 
\begin{eqnarray}
H_\alpha(w_\alpha) &:=& \frac14\sum_{n=3}^\infty
\frac{(-1)^n}{(2n-1)!}2^{2n-1}\alpha^{2n}(1-w_\alpha(x)^{2n-1}),
\\
&=& \sum_{n=3}^\infty
\frac{(-1)^n}{(2n-1)!}2^{2n-3}\alpha^{2n}(1-w_\alpha(x)^{2n-1}),
\nonumber
\\
J_\alpha(w_\alpha) &=& - \frac18\sum_{n=3}^\infty \frac{(-1)^n}{(2n)!}2^{2n}\alpha^{2n}
(1 - w_\alpha(x)^{2n})
\\
&=& - \sum_{n=3}^\infty \frac{(-1)^n}{(2n)!}2^{2n-3}\alpha^{2n}
(1 - w_\alpha(x)^{2n}).
\nonumber
\end{eqnarray}
We put 
\begin{eqnarray}
M_\alpha(w_\alpha)&:=& H_\alpha(w_\alpha) + J_\alpha(w_\alpha(x))
\\
&=& \sum_{n=3}^\infty(-1)^n2^{2n-3}
\left\{\frac{1}{(2n-1)!}(1-w_\alpha(x)^{2n-1}) - \frac{1}{(2n)!}(1-w_\alpha(x)^{2n})\right\}
\alpha^{2n}.
\nonumber
\end{eqnarray} 
By this and (2.3), for $0 \le x \le 1/2$, we have 
\begin{eqnarray}
w_\alpha'(x) = \sqrt{2\mu}\alpha^{-1}\sqrt{\frac{1}{p+1}\alpha^{p+1}(1-w_\alpha(x)^{p+1}) + \frac14\alpha^4
(1-w_\alpha(x)^4) + M_\alpha(w_\alpha)}.
\end{eqnarray} 
(i) Let $1 < p < 3$. Then by (3.7), we have 
\begin{eqnarray}
w_\alpha'(x) &=& \sqrt{2\mu\alpha^{-2}}\sqrt{\frac{\alpha^{p+1}}{p+1}}\sqrt{1-w_\alpha(x)^{p+1}}
\\
&&\times
\sqrt{ 1 + \frac{p+1}{4}\alpha^{3-p}\frac{1-w_\alpha(x)^4}{1-w_\alpha(x)^{p+1}} 
+ K(w_\alpha)\alpha^{5-p}},
\nonumber
\end{eqnarray}
where
\begin{eqnarray}
K(w_\alpha(x)):= - 2^3(p+1)\left\{\frac{1}{5!}\frac{1-w_\alpha(x)^5}{1-w_\alpha(x)^{p+1}} 
- \frac{1}{6!}\frac{1-w_\alpha(x)^6}{1-w_\alpha(x)^{p+1}} \right\}.
\end{eqnarray}
By (3.8) and Taylor expansion, we have 
\begin{eqnarray}
&&\sqrt{\frac{\mu}{2(p+1)}}\alpha^{(p-1)/2} = 
\\
&&= 
\int_0^{1/2} \frac{w_\alpha'(x)}{\sqrt{1-w_\alpha(x)^{p+1}}
\sqrt{ 1 + \frac{p+1}{4}\alpha^{3-p}\frac{1-w_\alpha(x)^4}{1-w_\alpha(x)^{p+1}} 
+ K(w_\alpha(x))\alpha^{5-p}}}dx
\nonumber
\\
&&=\int_0^1 \frac{1}{\sqrt{1-s^{p+1}}}
\left\{1 -\frac{p+1}{8}\alpha^{3-p}\frac{1-s^4}{1-s^{p+1}} 
-\frac12 K(s)\alpha^{5-p} + O(\alpha^{5-p})\right\}ds.
\nonumber
\end{eqnarray}
This implies from (1.7), (1.38) and (1.39) that 
\begin{eqnarray}
\sqrt{\mu} &=& \sqrt{2(p+1)}\alpha^{-(p-1)/2}\left\{C_{0,p} + L_1 \alpha^{3-p} 
+ L_2\alpha^{5-p} + O(\alpha^{7-p})\right\}.
\end{eqnarray}
This implies (1.41).

\noindent
(ii) Let $p = 3$. Then by (3.7), we have 
\begin{eqnarray}
w_\alpha'(x) &=& \sqrt{2\mu}\alpha^{-1}\sqrt{\frac12\alpha^4(1-w_\alpha(x)^4) 
+ M_\alpha(w_\alpha(x))}
\\
&=& \sqrt{\mu}\alpha\sqrt{1-w_\alpha(x)^4}\sqrt{1 + 
2\alpha^{-4}\frac{M_\alpha(w_\alpha(x))}{1-w_\alpha(x)^4}}.
\nonumber
\end{eqnarray}
This along with Taylor expansion implies that 
\begin{eqnarray}
\frac12\sqrt{\mu} &=& \alpha^{-1}\int_0^{1/2}\frac{w_\alpha'(x)}
{\sqrt{1-w_\alpha(x)^4}\sqrt{1 + 2\alpha^{-4}
\frac{M_\alpha(w_\alpha(x))}{1-w_\alpha(x)^4}}}dx
\nonumber
\\
&=& \alpha^{-1}\int_0^1 \frac{1}{\sqrt{1-s^4}}
\left\{1 - \alpha^{-4}\frac{M_\alpha(s)}{1-s^4} + O(\alpha^4)\right\}ds.
\nonumber
\end{eqnarray}
By this, we obtain 
\begin{eqnarray}
\sqrt{\mu} &=& 2\alpha^{-1}\int_0^1  \frac{1}{\sqrt{1-s^4}}
\left\{1 + 8\alpha^2\left(\frac{1}{5!}\frac{1-s^5}{1-s^4} - \frac{1}{6!}\frac{1-s^6}{1-s^4}\right)
+ O(\alpha^4)\right\}ds 
\nonumber
\\
&=& 2\alpha^{-1}\left\{C_{0,3} + \frac12H_3\alpha^{2} + O(\alpha^4)\right\}.
\nonumber
\end{eqnarray}
This implies (1.42). 

\noindent
(iii) 
Let $3 < p\le 5$. Then by (3.7), we have 
\begin{eqnarray}
\frac12\sqrt{2\mu} &=& 2\alpha\int_0^{1/2} \frac{w_\alpha'(x)}
{\alpha^2\sqrt{1-w_\alpha(x)^4}\sqrt{1 + \frac{4}{p+1}\alpha^{p-3}
\frac{1-w_\alpha(x)^{p+1}}{1-w_\alpha(x)^4} + Q_\alpha(w_\alpha(x))}},
\end{eqnarray}
where
\begin{eqnarray}
Q_\alpha(w_\alpha):= 4\alpha^{-4}\frac{M_\alpha(w_\alpha)}
{1-w_\alpha(x)^4}.
\end{eqnarray}
By this and Taylor expansion, we have 
\begin{eqnarray}
\sqrt{\frac{\mu}{2}} &=& 2\alpha^{-1}\int_0^1 \frac{1}{\sqrt{1-s^4}}
\\
&&\qquad \quad\times
\left\{1 - \frac{2}{p+1}\alpha^{p-3}\frac{1-s^{p+1}}{1-s^4} +\alpha^2\frac{2^4}{5!}
\frac{1-s^5}{1-s^4} -\alpha^2\frac{2^4}{6!}\frac{1-s^6}{1-s^4} 
+ O(\alpha^4) \right\}ds
\nonumber
\\
&=& 2\alpha^{-1}\left\{C_{0,3} + H_2\alpha^{p-3} + H_3\alpha^2 + O(\alpha^4)\right\}.
\nonumber
\end{eqnarray}
This implies (1.43). 

\noindent
(iv) Assume that $p > 5$. Then by (3.15), we have 
\begin{eqnarray}
\sqrt{\frac{\mu}{2}} &=& 2\alpha^{-1}\left\{C_{0,3} + H_3\alpha^2 + o(\alpha^2)\right\}.
\end{eqnarray}
This implies (1.44). Thus the proof is complete.

\section{Proof of Theorem 8} In this section, we assume that $0 < \alpha \ll 1$.  By Taylor expansion, we have 
\begin{eqnarray}
v_\alpha(x)\sin^2v_\alpha(x) &=& \sum_{n=1}^\infty \frac{(-1)^{n-1}2^{2n-1}}{(2n)!}
v_\alpha(x)^{2n+1}
\\
&=& v_\alpha(x)^3 - \frac13v_\alpha(x)^5 + \frac{2}{45}v_\alpha(x)^7 + O(v_\alpha(x)^9).
\nonumber
\end{eqnarray}

 \noindent
 (i) Let $1 < p < 3$. Then by (2.8), (4.1), Taylor expansion and 
 putting $v_\alpha = \theta = \alpha s$, we have 
 \begin{eqnarray}
 &&\Vert v_\alpha'\Vert_2^2 = 2\int_0^{1/2} v_\alpha'(x)v_\alpha'(x)dx 
 \\
 &&= 2\sqrt{2\mu}\int_0^{1/2}\sqrt{\frac{1}{p+1}(\alpha^{p+1}-v_\alpha(x)^{p+1}) 
 + \frac14(\alpha^4 - v_\alpha(x)^4)(1 + o(1))}
 v_\alpha'(x)dx
 \nonumber
 \\
 &&=  2\sqrt{2\mu}\int_0^{\alpha}\sqrt{\frac{1}{p+1}(\alpha^{p+1}-\theta^{p+1}) 
 + \frac14(\alpha^4 - \theta^4) (1 + o(1))}d\theta
 \nonumber
 \\
 &&= 2\sqrt{\frac{2\mu}{p+1}}\alpha^{(p+3)/2}\int_0^{1}\sqrt{1-s^{p+1}}
 \sqrt{1 + \frac{p+1}{4}\alpha^{3-p}\frac{1 - s^4}{1-s^{p+1}} (1 + o(1))}ds
 \nonumber
 \\
 &&= 2\sqrt{\frac{2\mu}{p+1}}\alpha^{(p+3)/2}\int_0^{1}\sqrt{1-s^{p+1}}
 \left\{1 +  \frac{p+1}{8}\alpha^{3-p}\frac{1 - s^4}{1-s^{p+1}} (1 + o(1))\right\}ds
 \nonumber
 \\
 &&= 2\sqrt{\frac{2}{p+1}}\sqrt{\mu}\alpha^{(p+3)/2}\left\{E_{0,p} + E_1\alpha^{3-p} + o(\alpha^{3-p})\right\}.
 \nonumber
 \end{eqnarray}
 By this and (3.11), we have 
 \begin{eqnarray}
 \Vert v_\alpha'\Vert_2^2  &=& 2\sqrt{\frac{2}{p+1}}\alpha^{(p+3)/2}
 \left\{E_{0,p} +E_1\alpha^{3-p} + o(\alpha^{3-p})\right\}
 \\
 &&\times \sqrt{2(p+1)}\alpha^{-(p-1)/2}\left\{C_{0,p} + L_1 \alpha^{3-p} 
+ L_2\alpha^{5-p} + O(\alpha^{7-p})\right\}
\nonumber
\\
&=& 4\alpha^2
\left\{E_{0,p}C_{0,p} + (E_{0,p}L_1+C_{0,p}E_1)\alpha^{3-p} + o(\alpha^{3-p})\right\}.
\nonumber
 \end{eqnarray}
 By this, (1.5) and Theorem 1.7 (i), we have 
 \begin{eqnarray}
 \lambda(\alpha) &=& 2(p+1)\alpha^{1-p}\left\{C_{0,p} + L_1\alpha^{3-p} 
+ L_2\alpha^{5-p} + O(\alpha^{7-p})\right\}^2
\\
&&\times
\left\{4b\alpha^2
\left\{E_{0,p}C_{0,p} + (E_{0,p}L_1+C_{0,p}E_1)\alpha^{3-p} + o(\alpha^{3-p})\right\} + 1\right\}.
\nonumber
 \end{eqnarray}
 (ii) Let $p = 3$. Then by (4.2) and putting $s = v_\alpha(x)/\alpha$, we have 
 \begin{eqnarray}
 \Vert v_\alpha'\Vert_2^2 &=& 2\sqrt{\mu}(1 + o(1))\int_0^{1/2} \sqrt{\alpha^4-v_\alpha(x)^4}
 v_\alpha'(x)dx 
 \\
 &=& 2\sqrt{\mu}(1 + o(1))\alpha^3\int_0^{1} \sqrt{1-s^4}ds
 \nonumber
 \\
 &=& 2\sqrt{\mu}\alpha^3E_{0,3}(1 + o(1)).
 \nonumber
 \end{eqnarray}
 By this, (1.5) and Theorem 1.7 (ii), we have 
 \begin{eqnarray}
 \Vert v_\alpha'\Vert_2^2 &=& 2\alpha^3E_{0,3}(1 + o(1))2\alpha^{-1}
 \left\{C_{0,3} + \frac12H_3\alpha^2 + O(\alpha^4)\right\}
 \\
 &=& 4\alpha^2E_{0,3}C_{0,3}(1 +o(1)).
 \nonumber
 \end{eqnarray}
 By this and Theorem 1.7 (ii), we have 
 \begin{eqnarray}
 \lambda(\alpha) &=& 4\alpha^{-2}(1 + 4bE_{0.3}C_{0,3}\alpha^2 + o(\alpha^2))
 \left\{C_{0,3} + \frac12H_3\alpha^2 + O(\alpha^4)\right\}^2.
 \end{eqnarray}
 
 \noindent
 We next consider the case $p > 3$. By (4.1), for $0 <  x < 1/2$, we have 
 \begin{eqnarray}
\frac12 v_\alpha'(x)^2 + \mu\left\{\frac14v_\alpha(x)^4 + \frac{1}{p+1}v_\alpha(x)^{p+1}
(1 + o(1))\right\} = \mu\left\{\frac14\alpha^4+ \frac{1}{p+1}\alpha^{p+1}
(1 + o(1))\right\}.
\end{eqnarray}
By this, for $0 \le x \le 1/2$, we have 
\begin{eqnarray}
v_\alpha'(x) &=& \sqrt{\frac{\mu}{2}}\sqrt{\alpha^4-v_\alpha(x)^4}
\sqrt{1 + \frac{4}{p+1}\frac{\alpha^{p+1} - v_\alpha(x)^{p+1}}{\alpha^4-v_\alpha(x)^4}(1 + o(1))}
\end{eqnarray}
By this, (3.16) and the same calculation as that of (4.2) and putting $v_\alpha(x) = \alpha s$, we have 
\begin{eqnarray}
\Vert v_\alpha'\Vert_2^2 &=& \sqrt{2\mu}\alpha^3
\int_0^1 \sqrt{1-s^4}\sqrt{1 + \frac{4}{p+1}\alpha^{p-3}\frac{1-s^{p+1}}{1-s^4}(1 + o(1))}ds
\\
&=& \sqrt{2\mu}\alpha^3\left\{E_{0,3} + E_5\alpha^{p-3}(1 + o(1))\right\}
\nonumber
\\
&=& 4\alpha^{2}\left\{C_{0,3} + H_2\alpha^{p-3} + H_3\alpha^2 + O(\alpha^4)\right\}
\left\{E_{0,p} + E_5\alpha^{p-3}(1 + o(1))\right\}.
\nonumber
\end{eqnarray}
(iii) Let $3 < p \le 5$. Then by (1.5), (3.17) and (4.10), we have 
\begin{eqnarray}
\lambda(\alpha) &=& 8\alpha^{-2}
\left\{C_{0,3} + H_2\alpha^{p-3} + H_3\alpha^2 + O(\alpha^4)\right\}^2
\\
&&\times
\left[4b\alpha^{2}\left\{C_{0,3} + H_2\alpha^{p-3} + H_3\alpha^2 + O(\alpha^4)\right\}
\left\{E_{0,3} + E_5\alpha^{p-3}(1 + o(1))\right\} + 1\right].
\nonumber
\end{eqnarray}
(iv) Let $p > 5$. Then by (1.27), (3.16) and  (4.10), we have 
\begin{eqnarray}
\lambda(\alpha) &=& 8\alpha^{-2}\left\{C_{0,3} + H_3\alpha^2 + o(\alpha^2)\right\}^2
\\
&&\times
\left[4b\alpha^{2}\left\{C_{0,3} + H_3\alpha^2 + o(\alpha^2)\right\}
\left\{E_{0,3} + E_5\alpha^{p-3}(1 + o(1))\right\} + 1\right].
\nonumber
\end{eqnarray}
Thus the proof of Theorem 1.8 is complete. \qed

\end{document}